\documentclass{article}
\usepackage{amssymb,amsfonts,amsmath,amsthm,enumerate}
\usepackage[numeric]{amsrefs}
\usepackage{mathtools}

\def \C{{\mathbb C}}

\def \f{\frac}

\def \G{\Gamma}

\def \Prim{\mathrm{Prim}}

\def \Z{{\mathbb Z}}
\def \l{\left}
\def \r{\right}

\newtheorem{theorem}{Theorem}

\numberwithin{equation}{section} \theoremstyle{definition}

\begin{document}

\title{Simple functional equations for generalized Selberg zeta functions with Tate motives}
\author{Shin-ya Koyama\footnote{Department of Biomedical Engineering, Toyo University,
2100 Kujirai, Kawagoe, Saitama, 350-8585, Japan.} \ \& Nobushige Kurokawa\footnote{Department of Mathematics, Tokyo Institute of Technology, 
Oh-okayama, Meguro-ku, Tokyo, 152-8551, Japan.}}

\maketitle

\begin{abstract}
We prove that for a compact locally symmetric Riemannian space $M$ of rank 1
there exist infinitely many automorphic Tate motives $f$ such that
the generalized Selberg zeta function $Z_{M(f)}(s)$ satisfies a simple functional equation
in the sense that it has no gamma factors.
\end{abstract}



\section*{Introduction}
Let
$$
Z_M(s)=\prod_{P\in\Prim(M)}\prod_{\lambda}\l(1-N(P)^{-s-\lambda}\r),
$$
be the Selberg zeta function of a compact locally symmetric Riemannian space $M$ of rank 1,
where $\Prim(M)$ denotes the set of primitive closed geodesics on $M$ with 
$
N(P)=\exp(\mathrm{length}(P)).
$
It is known that $Z_M(s)$ has analytic continuation to all $s\in\C$ with a functional equation
$$
Z_M(2\rho_M-s)=Z_M(s)S_M(s),
$$
where 
$$
S_M(s)=\exp\l(\mathrm{vol}(M)\int_0^{s-\rho_M}\mu_M(it)dt\r).
$$
Here $\mu_M(t)$ is the Plancherel measure and $2\rho_M,\ \mathrm{vol}(M)\in\Z_{>0}$.

Actually, the ``gamma factor'' $S_M(s)$ can be described explicitly in terms of the multiple sine function
(or the multiple gamma function) as in \cite{KK1}; see \S1 of the text.
In this paper we tacitly assume that the dimension $\mathrm{dim}(M)$ is even
since otherwise (that is the case of $M=\G\backslash SO(1,2n+1)/K$)
the factor $S_M(s)$ is a simple exponential function.

Now, we introduce the action of $f(x)\in\Z[x,x^{-1}]$ on $Z_M(s)$ as the
``(absolute) Tate motif'' (see Manin \cite{M}). 
More concretely, for a Laurent polynomial
$$
f(x)=\sum_{k\in\Z}a(k)x^k\in\Z[x,x^{-1}],
$$
we define
$$
Z_{M(f)}(s)=\prod_k Z_M(s-k)^{a(k)},
$$
which is considered as a generalized Selberg zeta function.

The purpose of this paper is to prove the following result.

\noindent{\bf Theorem}
Let $M$ be a compact locally symmetric Riemannian space of rank 1 as above.
Then $Z_{M(f)}(s)$ has simple functional equations
$$
Z_{M(f)}(D+2\rho_M-s)^C=Z_{M(f)}(s)
$$
for infinitely many automorphic $f$.
Here the automorphy of $f$ means that
$$
f(x^{-1})=Cx^{-D}f(x)
$$
with $C=\pm1$ and $D\in\Z$.

Our proof shows that the functional equation of $Z_{M(f)}(s)$ for an automorphic
$$
f(x)=\sum_{k\in\Z}a(k)x^k\in\Z[x,x^{-1}],
$$
is given as
$$
Z_{M(f)}(D+2\rho_M-s)^C=Z_{M(f)}(s)S_{M(f)}(s)
$$
where
$$
S_{M(f)}(s)=\prod_{k\in\Z} S_{M}(s-k)^{a(k)}
$$
in general: see Theorem 2 in \S2.
Hence the crucial point is to show that
$$
S_{M(f)}(s)=1
$$
for infinitely many $f$. Thus we do not need any ``gamma factors'' in such cases.
For more explicit construction of such $f$ see Theorem 3 in \S3.

We notice that the original Selberg zeta function $Z_M(s)$ was studied by Selberg \cite{S1, S2}
for a compact Riemann surface $M$ of genus $g\ge 2$.
It is written as
$$
Z_M(s)=\prod_{P\in\Prim(M)}\prod_{n=0}^\infty\l(1-N(P)^{-s-n}\r).
$$
The functional equation is
$$
Z_{M}(1-s)=Z_{M}(s)S_{M}(s)
$$
with
\begin{align*}
S_M(s)
&=\exp\l((4-4g)\int_0^{s-\f12}\pi t \tan(\pi t)dt\r)\\
&=(S_2(s)S_2(s+1))^{2-2g},
\end{align*}
where
$$
S_r(s)=S_r(s,(1,\cdots,1))
$$
is the multiple sine function in \cite{KK1}.
We refer to the paper \cite{KK2} concerning this case.

\section{The gamma factor}
We recall the explicit calculation of $S_M(s)$ by multiple sine functions and multiple gamma functions in \cite{KK1}
according to the classification of $M=\G\backslash G/K$.

\begin{theorem}[\cite{KK1}]
The factors $S_M(s)$ are written explicitly as follows.
\begin{enumerate}[\rm (1)]
\item 
$$
S_M(s)=\l\{\begin{array}{ll}
\l(S_{2n}(s)S_{2n}(s+1)\r)^{\mathrm{vol}(M)(-1)^{\f{\dim(M)}2}} & G=SO(1,2n)\\
\l(\prod\limits_{k=0}^nS_{2n}(s+k)^{\binom nk^2}\r)^{\mathrm{vol}(M)(-1)^{\f{\dim(M)}2}} & G=SU(1,n)\\
\l(\prod\limits_{k=0}^{2n-1}S_{4n}(s+k)^{\frac1{2n}\binom{2n}k\binom{2n}{k+1}}\r)^{\mathrm{vol}(M)}
& G=Sp(1,n)\\
(S_{16}(s)S_{16}(s+1)^{10}S_{16}(s+2)^{28}&\\
\quad\times S_{16}(s+3)^{28}S_{16}(s+4)^{10}S_{16}(s+5) )^{\mathrm{vol}(M)} & G=F_4.
\end{array}\r.
$$
\item
$$
S_M(s)=\f{\G_M(s)}{\G_M(2\rho_M-s)}
$$
with
$$
\G_M(s)=\l\{\begin{array}{ll}
\l(\G_{2n}(s)\G_{2n}(s+1)\r)^{\mathrm{vol}(M)(-1)^{\f{\dim(M)}2-1}} & G=SO(1,2n)\\
\l(\prod\limits_{k=0}^n\G_{2n}(s+k)^{{n \choose k}^2}\r)^{\mathrm{vol}(M)(-1)^{\f{\dim(M)}2-1}} & G=SU(1,n)\\
\l(\prod\limits_{k=0}^{2n-1}\G_{4n}(s+k)^{\frac1{2n}{2n \choose k}{2n \choose k+1}}\r)^{-\mathrm{vol}(M)}& G=Sp(1,n)\\
(\G_{16}(s)\G_{16}(s+1)^{10}\G_{16}(s+2)^{28} &\\
\quad\times\G_{16}(s+3)^{28}\G_{16}(s+4)^{10}\G_{16}(s+5))^{-\mathrm{vol}(M)} & G=F_4.
\end{array}\r.
$$
\end{enumerate}
\end{theorem}

We refer to \cite{KK1} for the proofs.
In the case of compact Riemann surfaces $M$ of genus $g\ge2$, we have
$$
M=\G\backslash SO(1,2)/K
$$
and $\mathrm{vol}(M)=2g-2$ with $\dim(M)=2$.

\section{Generalized functional equations}
We show the functional equation for $Z_{M(f)}(s)$ with automorphic Tate motif $f(x)\in\Z[x,x^{-1}]$
satisfying
$$
f(x^{-1})=Cx^{-D}f(x).
$$

\begin{theorem}
For each automorphic $f$ we have the functional equation
$$
Z_{M(f)}(D+2\rho_M-s)^C=Z_{M(f)}(s)S_{M(f)}(s).
$$
\end{theorem}

{\it Proof.}
First we remark that the following equivalence is easily seen:
$$
f(x^{-1})=Cx^{-D}f(x)\quad\Longleftrightarrow\quad a(D-k)=C a(k)\quad(\forall k\in\Z).
$$
Then we see that
\begin{align*}
Z_{M(f)}(D+2\rho_M-s)^C
&=\prod_k Z_{M}(D+2\rho_M-s-k)^{Ca(k)}\\
&=\prod_k Z_{M}(2\rho_M-s+(D-k))^{a(D-k)}\\
&=\prod_k Z_{M}(2\rho_M-(s-k))^{a(k)},
\end{align*}
where we replaced $k$ by $D-k$. Thus the functional equation for $Z_M(s)$ gives
\begin{align*}
Z_{M(f)}(D+2\rho_M-s)^C
&=\prod_k (Z_{M}(s-k)S_{M}(s-k))^{a(k)}\\
&=Z_{M(f)}(s)S_{M(f)}(s).
\end{align*}
\hfill\qed

\section{Vanishment of gamma factors}
Now we prove the main result.
\begin{theorem}
Let $f(x)\in(x-1)^{\dim(M)}\Z[x,x^{-1}]$ be automorphic. Then it holds that
$$
S_{M(f)}(s)=1.
$$
Namely, we have 
$$
Z_{M(f)}(D+2\rho_M-s)^C=Z_{M(f)}(s)
$$
for such $f$.
\end{theorem}

{\it Proof.}
We denote by $S_r^f(s)$ the action of $f(x)\in\Z[x,x^{-1}]$ on $S_r(x)$. 
We remark that
$$
S_r^{fg}=(S_r^f)^g.
$$
In fact for
\begin{align*}
f(x)&=\sum_k a(k)x^k,\\
g(x)&=\sum_l b(l)x^l
\end{align*}
it holds that
$$
f(x)g(x)=\sum_{k,l}a(k)b(l)x^{k+l}
$$
and
$$
S_r^{fg}(s)=\prod_{k,l}S_r(s-(k+l))^{a(k)b(l)}.
$$
On the other hand it follows that
$$
S_r^f(s)=\prod_k S_r(s-k)^{a(k)}
$$
and that
\begin{align*}
(S_r^f)^g(s)
&=\prod_l\l(\prod_k S_r((s-l)-k)^{a(k)}\r)^{b(l)}\\
&=\prod_{k,l}S_r(s-(k+l))^{a(k)b(l)}.
\end{align*}
Thus $S_r^{fg}=(S_r^f)^g$.

Let $r=\dim(M)$ and put
$$
f_r(x)=(1-x^{-1})^r.
$$
We show that
$$
S_r^{f_r}(s)=-1.
$$
The case $r=1$ is easy:
$$
S_1^{f_1}(s)
=\f{S_1(s)}{S_1(s+1)}=\f{2\sin(\pi s)}{2\sin(\pi(s+1))}
=-1.
$$
In general, the relation (see \cite{KK1})
$$
S_r^{(1-x^{-1})}(s)
=\f{S_r(s)}{S_r(s+1)}=S_{r-1}(s)
$$
gives
$$
S_r^{f_r}(s)=S_{r-1}^{f_{r-1}}(s)=\cdots=S_1^{f_1}(s)=-1
$$

Let
$$
f(x)\in(x-1)^{r}\Z[x,x^{-1}]
$$
and
$$
f(x)=f_r(x)g(x).
$$
Then we have
$$
S_r^f(s)=\pm1.
$$
Actually,
$$
S_r^f(s)=S_r^{f_{r}g}(s)=(S_r^{f_r})^g(s)=(-1)^{g(1)}=\pm1.
$$
Thus by the explicit formula for $S_M(s)$ in \S1, we see that
$$
S_{M(f)}(s)=1.
$$
\hfill\qed

{\bf Remark.}
The above arguments give the implication 
\begin{align*}
f(1)=f'(1)=\cdots=f^{(r-1)}(1)=0&\\
\Longrightarrow\quad& S_r^f(s)=\pm1.
\end{align*}
We refer to \cite{KT, KT2} concerning the converse and generalizations.


\begin{bibdiv} \begin{biblist}

\bib{KK1}{article}{
   author={S. Koyama},
   author={N. Kurokawa},
   title={Multiple sine functions},
   journal={Forum Math.},
   volume={15},
   date={2003},
   pages={839--876},
}
\bib{KK2}{article}{
   author={S. Koyama},
   author={N. Kurokawa},
   title={Functional equations for Selberg zeta functions with Tate motives},
   journal={(preprint)},
   date={2020},
}
\bib{KT}{article}{
   author={Kurokawa, N.},
   author={H. Tanaka},
   title={Absolute zeta functions and the automorphy},
   journal={Kodai Math. J.},
   volume={40},
   date={2017},
   pages={584-614},
}
\bib{KT2}{article}{
   author={Kurokawa, N.},
   author={H. Tanaka},
   title={Reductions of multiple sine functions and multiple gamma functions},
   journal={(in preparation)},
}
\bib{M}{article}{
   author={Y. Manin},
   title={Lectures on zeta functions and motives (according to Deninger and Kurokawa)},
   journal={Ast\'erisque},
   volume={228},
   date={1995},
   pages={121-163},
}
\bib{S1}{article}{
   author={Selberg, A.},
   title={Harmonic analysis and discontinuous groups in weakly symmetric Riemannian spaces with applications to Dirichlet series},
   journal={J. Indian Math. Soc.},
   volume={20},
   date={1956},
   pages={47-87},
}
\bib{S2}{inproceedings}{
   author={Selberg, A.},
   title={G\"ottingen lectures},
   publisher={Springer Verlag}
   booktitle={Collected Works, Vol. I}
   date={1989},
   pages={626-674},
}
\end{biblist} \end{bibdiv}
\end{document}